\documentclass[12pt,leqno]{article}
\usepackage{amsfonts}
\usepackage{amsmath, amsthm, amsfonts, amssymb, color}
\usepackage{mathrsfs}
\usepackage{url}
\usepackage{color}
\theoremstyle{definition}

\newcommand{\scr}[1]{\mathscr #1}
\definecolor{wco}{rgb}{0.5,0.2,0.3}

\numberwithin{equation}{section} \theoremstyle{remark}

\newcommand{\ua}{\uparrow}

\title{{\bf  Convergence in Wasserstein Distance for Empirical Measures of Semilinear SPDEs }\footnote{Supported in
 part by   NNSFC (11831014, 11921001) and   the DFG through the CRC 1283.} }
\author{
{\bf    Feng-Yu Wang$^{a),b)}$    }\\
\footnotesize{$^{a)}$ Center for Applied Mathematics, Tianjin University, Tianjin 300072, China }\\
 \footnotesize{ $^{b)}$ Department of Mathematics,
Swansea University,
Bay Campus,
Swansea,
SA1 8EN, United Kingdom}  }

\begin{document}
\allowdisplaybreaks
\def\R{\mathbb R}  \def\ff{\frac} \def\ss{\sqrt} \def\B{\mathbf
B}\def\TO{\mathbb T}
\def\I{\mathbb I_{\pp M}}\def\p<{\preceq}
\def\N{\mathbb N} \def\kk{\kappa} \def\m{{\bf m}}
\def\ee{\varepsilon}\def\ddd{D^*}
\def\dd{\delta} \def\DD{\Delta} \def\vv{\varepsilon} \def\rr{\rho}
\def\<{\langle} \def\>{\rangle} \def\GG{\Gamma} \def\gg{\gamma}
  \def\nn{\nabla} \def\pp{\partial} \def\E{\mathbb E}
\def\d{\text{\rm{d}}} \def\bb{\beta} \def\aa{\alpha} \def\D{\scr D}
  \def\si{\sigma} \def\ess{\text{\rm{ess}}}
\def\beg{\begin} \def\beq{\begin{equation}}  \def\F{\scr F}
\def\Ric{{\rm Ric}} \def\Hess{\text{\rm{Hess}}}
\def\e{\text{\rm{e}}} \def\ua{\underline a} \def\OO{\Omega}  \def\oo{\omega}
 \def\tt{\tilde}
\def\cut{\text{\rm{cut}}} \def\P{\mathbb P} \def\ifn{I_n(f^{\bigotimes n})}
\def\C{\scr C}      \def\aaa{\mathbf{r}}     \def\r{r}
\def\gap{\text{\rm{gap}}} \def\prr{\pi_{{\bf m},\varrho}}  \def\r{\mathbf r}
\def\Z{\mathbb Z} \def\vrr{\varrho} \def\ll{\lambda}
\def\L{\scr L}\def\Tt{\tt} \def\TT{\tt}\def\II{\mathbb I}
\def\i{{\rm in}}\def\Sect{{\rm Sect}}  \def\H{\mathbb H}
\def\M{\scr M}\def\Q{\mathbb Q} \def\texto{\text{o}} \def\LL{\Lambda}
\def\Rank{{\rm Rank}} \def\B{\scr B} \def\i{{\rm i}} \def\HR{\hat{\R}^d}
\def\to{\rightarrow}\def\l{\ell}\def\iint{\int}
\def\EE{\scr E}\def\Cut{{\rm Cut}}\def\W{\mathbb W}
\def\A{\scr A} \def\Lip{{\rm Lip}}\def\S{\mathbb S}
\def\BB{\scr B}\def\Ent{{\rm Ent}} \def\i{{\rm i}}\def\itparallel{{\it\parallel}}
\def\g{{\mathbf g}}\def\Sect{{\mathcal Sec}}\def\T{\mathcal T}\def\V{{\bf V}}
\def\PP{{\bf P}}\def\HL{{\bf L}}\def\Id{{\rm Id}}\def\f{{\bf f}}\def\cut{{\rm cut}}

\def\BL{\scr A}

\maketitle

\begin{abstract} The convergence rate in Wasserstein distance is estimated for the empirical measures of symmetric semilinear SPDEs. Unlike in the finite-dimensional case that the convergence is of  algebraic order in time, in the present situation the convergence is of log order with a power given by eigenvalues of the underlying linear operator.  
\end{abstract} \noindent
 \vskip 2cm

\section{Introduction}
As the continuous Markov process counterpart   of Wasserstein matching problem for i.i.d. samples studied in \cite{AMB, [25]} and references within, in \cite{W20, W20b, W20c, WZ20} we have estimated the convergence rate in Wasserstein distance for empirical measures
of symmetric diffusion processes.

Let $V\in C^2(M)$ for a $d$-dimensional compact connected Riemannian manifold $M$, let  $X_t$ be the diffusion process generated by $L:=\DD+\nn V$ on $M$ with reflecting boundary if exists, and let $\W_2$ be the $L^2$-Wasserstein distance induced by the Riemannian metric. According to \cite{WZ20},
the empirical measure $\mu_t:=\ff 1 t\int_0^t \dd_{X_s}\d s$ satisfies
$$\lim_{t\to\infty} t \E[\W_2(\mu_t,\mu)^2] =\sum_{i=1}^\infty \ff 2 {\dd_i},$$
where $\{\dd_i\}_{i\ge 1}$ are all non-trivial eigenvalues of   $-L$ in $L^2(\mu)$ counting multiplicities, with Neumann condition if the boundary exists. Since $\sum_{i=1}^\infty \ff 2 {\dd_i}<\infty$ if and only if $d\le 3$, so that when $t\to\infty$
 $$\E[\W_2(\mu_t,\mu)^2]\approx \ff 1 t,\ \ d\le 3,$$
 where we write $a(t) \approx b(t)$ for two positive functions $a$ and $b$ on $(0,\infty)$, if there exists a constant $C>1$ such that $C^{-1} a(t)\le b(t)\le C a(t)$ holds for large $t>0.$ 
 Moreover, we have proved in \cite{WZ20} that
 $$\E[\W_2(\mu_t,\mu)^2]\approx \beg{cases}\ff 1 t \log t,\ &\text{if}\ d=4,\\
 t^{-\ff 2{d-2}},\ &\text{if}\ d\ge 5.\end{cases} $$
 These results were  then extended in \cite{W20, W20b} for the empirical measure $\mu_t$ of conditional Dirichlet diffusion processes not  reaching the boundary before time $t$, and in \cite{W20c} for diffusion processes on non-compact complete Riemannian manifolds.

  In this paper, we
 investigate the problem for semilinear SPDEs, whose solutions provide a fundamental class of infinite-dimensional diffusion processes, see \cite{DZ1,DZ2} for details. It turns out that for this kind of infinite-dimensional processes the convergence of empirical measures becomes log order with a power determined by eigenvalues of the underlying linear operator.

Consider the following SDE on a separable Hilbert space $\H$:
\beq\label{E1} \d X_t= \big\{\nn V(X_t)-A X_t \big\}\d t + \ss 2\, \d W_t,\end{equation}
where $W_t$ is the cylindrical Brownian motion on $\H$, i.e. 
$$W_t= \sum_{i=1}^\infty B_t^i e_i,\ \ t\ge 0$$ for an orthonormal basis $\{e_i\}_{i\ge 1}$ of $\H$ and a sequence of independent one-dimensional Brownian motions $\{B_t^i\}_{i\ge 1},$ 
 $(A,\D(A))$ is a positive definite self-adjoint operator  and $V\in C^1(\H)$ satisfying the following assumption.
 
\beg{enumerate} \item[$(H_1)$] $A$ has discrete spectrum with eigenvalues $\{\ll_i>0\}_{i\ge 1}$ listed in the increasing order counting multiplicities  satisfying
$\sum_{i=1}^d \ll_i^{-\dd}<\infty$ for some constant $\dd\in (0,1)$, and $V\in C^1(\H)$, $\nn V$ is Lipschitz continuous in $\H$ such that \beq\label{VV} \<\nn V (x)- \nn V(y), x-y\>\le (K+\ll_1)|x-y|^2,\ \ x,y\in \H\end{equation}
holds for some constant $K\in\R$. Moreover,
$Z_V:= \mu_0(\e^V)<\infty,$  where
  $\mu_0$ is the centered Gaussian measure on $\H$ with covariance operator $A^{-1}$. \end{enumerate}
Under this condition,  for any $\F_0$-measurable random variable $X_0$ on $\H$, \eqref{E1} has a unique mild solution, and there exists an increasing function $\psi: [0,\infty)\to (0,\infty)$ such that
\beq\label{A00} \E[|X_t|^2]\le  \psi(t)\big(1+ \E[|X_0|^2]\big),\ \ t\ge 0,\end{equation}
see for instance \cite[Theorem 3.1.1]{W13}, or the earlier monographs
 \cite{DZ1,DZ2}.

 Let  $P_t$ be the associated Markov semigroup, i.e.
$$P_tf(x):= \E^x [f(X_t)],\ \ t\ge 0, f\in \B_b(\H),\ \ x\in \H,$$
where $\B_b(\H)$ is the class of all bounded measurable functions on $\H$, and $\E^x$ is   the expectation   for the solution $X_t$ of \eqref{E1} with $X_0=x$. In general, for a probability measure $\nu$ on $\H$, let $\E^\nu$ be the expectation for $X_t$ with initial distribution $\nu$.

By $(H_1),$ we define the probability measure
  $$\mu(\d x):= Z_V^{-1} \e^{V(x)}\mu_0(\d x).$$
Then  $P_t$   is symmetric in $L^2(\mu)$.  For any $p\ge 1$, the $L^p$-Wasserstein distance is given by
$$\W_p(\mu_1,\mu_2):= \inf_{\pi\in \C(\mu_1,\mu_2) }\bigg(\int_{\H\times\H} |x-y|^p\pi(\d x,\d y)\bigg)^{\ff 1 p},\ \ \mu_1,\mu_2\in \scr P(\H),$$
where $\scr P(\H)$ is the set of all probability measures on $\H$ and $\C(\mu_1,\mu_2)$ is the class of all couplings of $\mu_1$ and $\mu_2.$ 

In the following two sections, we investigate the upper bound and lower bound estimates on $\W_p(\mu_t,\mu)$ for the empirical measures
$$\mu_t:=\ff 1 t \int_0^t \dd_{X_s}\d s,\ \ t>0$$
of solutions to \eqref{E1}, 
where $\dd_x$ stands for the Dirac measure at point $x$.  Concrete examples are given to illustrate the resulting estimates, which show that in the present setting the convergence rate is of log order in $t$ with a power given by the growth of $\ll_i$ as $i\to\infty$. In particular, when $|V(x)|\le c(1+|x|)$ for some constant $c>0$ and all $x\in \H$, and $\ll_i \approx i^p$ for some $p>1$ and large $i$, Example 2.1 and Example 3.1 below imply
$$c_1 (\log t)^{1-p\land 3} \le \E^\mu\big[\W_2(\mu_t,\mu)^2\big] \le c_2 (\log t)^{\ff 1 p -1}$$ 
for some constants $c_1,c_2>0$ and large $t>0.$ 

To conclude this section, we compare the present study with the corresponding ones in \cite{W20}-\cite{WZ20} where finite-dimensional diffusion processes are investigated. 
Due to the lack of sharp estimates on heat kernel and eigenvalues for the generator,  new techniques have been developed for the present infinite-dimensional setting.  In particular, for the upper bound estimate we apply the dimension-free Harnack inequality established by the author (Section 2), while for the lower bound estimate we present a general result on the Wasserstein distance for discrete measures which applies to infinite-dimensions (Section 3). To derive a reasonable convergence rate,    optimization methods are applied for Wasserstein distances of regularized  and discretized approximations of empirical measures. 

\section{Upper bound estimate}

We first observe that $(H_1)$ implies   the following dimension-free Harnack inequality:
\beq\label{HN} (P_tf(x))^p\le (P_t f^p(y)) \exp\Big[\ff{pK|x-y|^2}{2(p-1) (1-\e^{-2Kt})}\Big],\ \ t>0, x,y\in \H, f\in \B^+(\H),\end{equation}
where $\B^+(\H)$ is the class of all nonnegative measurable  functions  on $\H$.   Indeed, by $(H_1)$, the operator $(\ll_1-A,\D(A))$ satisfies {\bf (A3.3)} in \cite{W13}, while $b(s,x):= \nn V(x) -\ll_1x$ and $S(t):= \e^{(\ll_1 -A)t}$ satisfy {\bf (A3.1)}  and {\bf (A3.2)}   in \cite{W13}. So, \eqref{HN} follows from  \cite[Theorem 3.2.1]{W13}.

Next, according to \cite[Theorem 1.4.1(6)]{W13}, \eqref{HN} implies that $P_t$ has a (symmetric) heat kernel $p_t(x,y)$ with respect to $\mu$ such that
\beq\label{AK} \beg{split} & \mu\big(p_t(x,\cdot)^{\ff p{p-1}}\big)^{p-1}=\sup_{\mu(|f|^p)\le 1} (P_tf(x))^p \\
&\le \bigg(\int_\H \e^{-\ff{pK|x-y|^2}{(p-1)(1-\e^{-2Kt})}}\mu(\d y)\bigg)^{-1},\ \ x\in \H, t>0, p>1.\end{split}\end{equation} 
In particular, by taking $p=2$ we obtain
 \beq\label{A0} p_{2t}(x,x)\le c(t,x):=   \bigg(\int_\H \e^{-\ff{2K|x-y|^2}{1-\e^{-2Kt}}}\mu(\d y)\bigg)^{-1}<\infty,\ \ t>0, x\in \H.\end{equation}  We assume that for any $t>0$,
\beq\label{AA}\beg{split} & \aa(t):= \E^\mu\big[|X_0-X_t|^2\big]= \int_{\H\times\H} |x-y|^2 p_t(x,y) \mu(\d x)\mu(\d y)<\infty,\\
& \bb(t):= \int_{\H} p_{2t}(x,x)\mu(\d x)= \int_{\H\times\H} p_t(x,y)^2 \mu(\d x)\mu(\d y)<\infty,\ \ t>0.\end{split}\end{equation}
In particular, $\bb(t)<\infty$ implies the uniform integrability of $P_t$ in $L^2(\mu)$, so that by \cite[Lemma 3.1]{GW},  $P_t$ is compact in $L^2(\mu)$ and the generator $L$ has discrete spectrum.
Since the associated Dirichlet form is irreducible, this implies that $L$ has a spectral gap $\ll_0>0$, such that 
\beq\label{*EXP} \mu(|P_tf-\mu(f)|^2) \le \e^{-2\ll_0 t} \mu(|f-\mu(f)|^2),\ \ t\ge 0, f\in L^2(\mu).\end{equation}
  In the following theorem,  we use $\aa$ and $\bb$ to estimate  the convergence rate of  $\E [\W_2(\mu_t,\mu)^2]$ as $t\to\infty$.

\beg{thm}\label{T1.1} Assume $(H_1)$ and  $\eqref{AA}$, and let $c(t,x)$ be in $\eqref{A0}$. We have
\beq\label{A-1} \E^\mu \big[\W_2(\mu_t,\mu)^2\big]\le   \inf_{\vv\in (0,1)} \Big\{\ff{16 \bb(\vv)}{\ll_0^2 t}  + 2 \aa(\vv)\Big\}=:\xi_t,\ \ t>0.\end{equation}
Consequently, for any $x\in \H$,
\beq\label{A2} \big(\E^x [\W_2(\mu_t,\mu)]\big)^2\le \inf_{r>0} \Big\{ \ff {8r} t \sup_{s\ge 0}  \E^x |X_s|^2 + 2 c(r,x) \xi_t \Big\},\ \ t>0.\end{equation}
\end{thm}
\beg{proof}
(a) We will use the following inequality due to \cite[Theorem 2]{Ledoux}:
\beq\label{Ledoux} \W_2(f\mu,\mu)^2\le 4 \mu(|\nn (-L)^{-1} (f-1)|^2),\ \ f\ge 0, \mu(f)=1.\end{equation} This estimate was 
  proved using the Kantorovich dual formula and the Hamilton-Jacobi equations, see  \cite{AMB} for an alternative estimate.

To apply \eqref{Ledoux}, we   consider    the modified empirical measures
\beq\label{B5} \mu_{\vv,t}:=\mu_t P_\vv = f_{\vv,t}\mu,\ \ \vv>0, t>0,\end{equation}
where
\beq\label{B6}  f_{\vv,t}:=\ff 1 t\int_0^t p_\vv(X_s,\cdot)\d s.\end{equation}
Noting that
$$P_s \{p_\vv(x,\cdot)\}(y)= p_{s+\vv}(x,y),\ \ x,y\in \H, s\ge 0,$$
by the spectral representation we obtain
\beq\label{BK1} \beg{split} &\mu(|\nn (-L)^{-1} (f_{\vv,t}-1)|^2) = \int_0^\infty \mu(|P_{s/2}(f_{\vv,t}-1)|^2)\d s\\
 & =\int_0^\infty \d s\int_\H \bigg(\ff 1 t\int_0^t  \big(p_{\vv+s/2}(X_u,\cdot)-1\big) \d u \bigg)^2 \d\mu\\
&= \ff 2 {t^2}\int_0^\infty \d s\int_0^t \d s_1 \int_{s_1}^t \mu \big(\{p_{\vv+s/2}(X_{s_1} ,\cdot)-1\}\cdot \{p_{\vv+s/2}(X_{s_2} ,\cdot)-1\}\big)\d s_2\\
 &=  \ff 2 {t^2}\int_0^\infty \d s\int_0^t \d s_1 \int_{s_1}^t  \big\{ p_{2\vv+s}(X_{s_1},X_{s_2})-1\big\}\d s_2.
 \end{split}\end{equation} 
 Next, by \eqref{*EXP} we have 
 \beq\label{BK2}  p_{r+s}(x ,x)-1=\mu\big(|P_{\ff s 2} \{p_{\ff r 2}(x,\cdot)\}-1|^2\big) \le \e^{-\ll_0 s} \big\{p_r(x,x)-1\big\},\ \ s,r>0.\end{equation} 
 Combining this with the Markov property we derive 
   \beg{align*} &\E^\mu  \big\{ p_{2\vv+s}(X_{s_1},X_{s_2})-1\big\}=\int_\H P_{s_2-s_1} \big\{p_{2\vv+s} (x,\cdot)-1\big\}(x) \mu(\d x)\\
 &= \int_\H   \big\{p_{2\vv+s+ s_2-s_1 } (x,x)-1\big\}  \mu(\d x)\le \e^{-\ll_0(s+s_2-s_1)} \bb(\vv).\end{align*}
 Therefore, \eqref{Ledoux} for $f:= f_{\vv,t}$ and \eqref{BK1} imply 
 \beq\label{A3}\beg{split} & \E^\mu \big[\W_2(\mu_{\vv,t},\mu)^2\big]\le \ff{8 \bb(\vv)}{t^2}\int_0^\infty\d s    \int_0^t \d s_1 \int_{s_1}^t \e^{-\ll_0(s+ s_2-s_1)}\d s_2 \\
 &= \ff{8 \bb(\vv)}{t^2}\int_0^\infty\d s    \int_0^t \ff{\e^{-\ll_0 s}-\e^{-\ll_0 (s+t-s_1)}}{\ll_0} \d s_1\le  \ff{8 \bb(\vv)}{t^2}\int_0^\infty\d s    \int_0^t \ff{\e^{-\ll_0 s}}{\ll_0} \d s_1\\
 &=  \ff{8 \bb(\vv)}{t\ll_0 }\int_0^\infty\e^{-\ll_0 s} \d s   
 = \ff {8 \bb(\vv)} {t\ll_0^2},\ \ t,\vv>0.\end{split}\end{equation}

 On the other hand, by  Jensen's inequality and that $ \dd_{x}P_\vv =\L_{X_\vv}$ for $X_0=x$, we obtain
 $$\W_2(\mu_{\vv,t},\mu_t)^2\le \bigg(\ff 1 t \int_0^t \W_2(\dd_{X_s}, \dd_{X_s}P_\vv)\d s\bigg)^2
 \le \ff 1 t \int_0^t  \big\{\E^x\big[ |x- X_\vv|^2 \big] \big\}\big|_{x=X_s}\d s.$$
 Since $\L_{X_s}=\mu$ for $\L_{X_0}=\mu$, this implies
 $$\E^\mu \big[\W_2(\mu_{\vv,t},\mu_t)^2\big]\le   \E^\mu\big[|X_\vv-X_0|^2\big]=\aa(\vv).$$
 Combining  with \eqref{A3}, we derive
 $$\E^\mu \big[\W_2(\mu_t,\mu)^2\big] \le 2 \E^\mu [\W_2(\mu_{\vv, t},\mu)^2]+ 2 \E^\mu [\W_2(\mu_{\vv,t},\mu)^2]\le \ff {16 \bb(\vv)} {t\ll_0} + 2 \aa(\vv),\ \ \vv\in (0,1).$$
Therefore, \eqref{A-1} holds.

(b) For any  $r>0$,  let  
$$ \mu_t^{(r)}:= \ff 1 t \int_r^{t+r} \dd_{X_s}\d s,\ \ t>0.$$
By the Markov property, Schwarz  inequality, \eqref{A0} and  \eqref{A-1}, we obtain
\beq\label{A4} \beg{split} &\big(\E^x [\W_2(\mu_t^{(r)},\mu)]\big)^2= \bigg(\int_\H \E^y [\W_2(\mu_t,\mu)] p_r(x,y)\mu(\d y)\bigg)^2\\
&\le p_{2r}(x,x) \int_\H \E^y \big[\W_2(\mu_t,\mu)^2\big] \mu(\d y) \\
&\le  c(x,r) \E^\mu \big[\W_2(\mu_t,\mu)^2\big] \le c(x,r) \xi_t,\ \ t>0.\end{split}\end{equation}
On the other hand, it is easy to see that
$$\pi_t :=  \ff 1 t\int_0^{r\land t}  \dd_{(X_s, X_{s+\ff{tr}{t\land r}})} \d s +\ff 1 t \int_{r\land t}^{t} \dd_{(X_s,X_s)}\d s\in \C(\mu_t,\mu_t^{(r)}), $$
so that
\beg{align*} &\E^x\big[\W_2(\mu_t,\mu_t^{(r)})^2 \big]\le \E^x \int_{\R^d\times\R^d}|y-z|^2 \pi_t(\d y,\d z)\\
& =\ff{1}{t} \int_0^{t\land r} \E^x |X_s-X_{r+\ff{ts}{r\land t}}|^2\d s\le \ff{4 r}t \sup_{s\ge 0} \E^x|X_s|^2.\end{align*}
This together with \eqref{A4} and the triangle inequality for $\W_2$, we prove \eqref{A2}.
 \end{proof}

Since the heat kernel $p_t(x,y)$ is usually unknown, the estimate presented in   Theorem \ref{T1.1}  is not explicit.  To derive explicit estimates,  we make the following assumption.

\beg{enumerate} \item[$(H_2)$]  There exists  an increasing function $\gg: (0,\infty)\to [0,\infty)$ such that
$$ |V(x)|\le \ff 1 2 \big(\gg(\vv^{-1})+  \vv   |x|^2\big),\ \ x\in \H,\vv>0.$$
\item[$(H_3)$] There exist   constants $c>0$ and   $\theta\in [0, \ll_1)$
$$ \ |\nn V(x)|\le c+ \theta |x|,\ \ x\in \H.$$
\end{enumerate}

\beg{cor}\label{C1.2} Assume $(H_1)$ and $(H_2)$.  Then: \beg{enumerate}
\item[$(1)$] There exists a constant $c_0>0$ such that
\beq\label{B1} \E^\mu\big[\W_2(\mu_t,\mu)^2\big]\le c_0 \inf_{\vv\in (0,1)}    \bigg( \ff 1 t \e^{k\vv^{-1}+\gg(k\vv^{-1})}  +\sum_{i=1}^\infty
\ff{1-\e^{-2\ll_i\vv}}{\ll_i} \bigg)=:\eta_t,\ \ t>0.\end{equation}
\item[$(2)$] If  $(H_3)$ holds,  then for any $k>K^+$ there exists a constant $c(k)>0$ such that
\beq\label{B2} \big(\E^x[\W_2(\mu_t,\mu)]\big)^2\le c(k)  \e^{k|x|^2} \eta_t,\ \ t\ge 1. \end{equation}\end{enumerate}
\end{cor}

To prove this result, we need the following two lemmas. 

\beg{lem} \label{L1}  Assume $(H_1)$ and $(H_3)$. There exists a constant $k>0$ such that
$$\sup_{t\ge 0} \E^x [|X_t|^2]\le k(1+|x|^2),\ \ x\in \H.$$
\end{lem}

\beg{proof} For $X_0=x$ we have
$$X_t= \e^{-A t} x +\int_0^t \e^{-A(t-s)}\nn V(X_s)\d s +\ss 2 \int_0^t \e^{-A(t-s)} \d W_s,\ \ t\ge 0.$$
By $(H_1)$ and $(H_3)$, we obtain
\beq\label{A*C} \beg{split} \E|X_t|^2 &\le (1+\vv^{-1}) \E^x\bigg|\e^{-At}x + \ss 2 \int_0^t \e^{A(t-s)}\d W_s\bigg|^2\\
&\qquad +(1+\vv) \E^x\bigg|\int_0^t\e^{-A(t-s)}\nn V(X_s)\d s\bigg|^2\\
&\le 2 (1+\vv^{-1}) \bigg(\big|\e^{-At}x\big|^2 + 2 \E^x \bigg|\int_0^t \e^{-A(t-s)}\d W_s\bigg|^2\bigg) \\
&\qquad + (1+\vv) \bigg(\int_0^t \big|\e^{-A(t-s)} \nn V(X_s) \big|\d s\bigg)^2\\
&\le 2 (1+\vv^{-1}) \bigg(\e^{-2\ll_1 t}|x|^2 + 2 \sum_{i=1}^\infty \int_0^t \e^{-2\ll_i(t-s)}\d s\bigg)\\
&\qquad + (1+\vv) \E^x\bigg(\int_0^t \e^{- \ll_1(t-s)} \big(c+\theta|X_s|\big)\d s\bigg)^2,\ \ \vv>0, t\ge 0.\end{split}\end{equation}
Noting that $(H_1)$ implies 
$$ \sum_{i=1}^\infty \int_0^t \e^{- 2\ll_i(t-s)}\d s\le \sum_{i=1}^\infty\ff 1 {2\ll_i}<\infty,$$   and that by Jensen's inequality
\beg{align*} &\E^x\bigg(\int_0^t \e^{- \ll_1(t-s)} \big(c+\theta|X_s|\big)\d s\bigg)^2\le \ff 1 { \ll_1}\E^x \int_0^t \e^{- \ll_1(t-s)} \big(c+\theta|X_s|\big)^2\d s\\
&\le \ff {(1+\vv^{-1})c^2} { \ll_1}+ \ff {(1 +\vv)\theta^2 }{ \ll_1} \int_0^t \e^{- \ll_1(t-s)} \E^x[|X_s|^2]\d s,\end{align*}
for any $\vv>0$ we find a constant $C(\vv)>0$ such that \eqref{A*C} yields
$$\E|X_t|^2
  \le C(\vv) (1+|x|^2)  +  \ff{(1+\vv)^2\theta^2}{2\ll_1} \int_0^t \e^{-\ll_1 (t-s)}  \E[|X_s|^2]\d s,\ \ \vv>0, t\ge 0.$$ 
Since $\theta< \ll_1$, we may take  $\vv>0$ such that $\ll_\vv:=\ff{(1+\vv)^2\theta^2}{ \ll_1}< \ll_1$, so that
$$\e^{ \ll_1 t} \E^x[|X_t|^2] \le C(\vv) (1+|x|^2) \e^{\ll_1 t} +\ll_\vv\int_0^t  \e^{ \ll_1 s} \E |X_s|^2 \d s,\ \ t\ge 0.$$
By Gronwall's lemma,  we obtain
\beg{align*} \e^{ \ll_1 t} \E^x[|X_t|^2] &\le C(\vv) (1+|x|^2) \e^{ \ll_1 t} +\ll_\vv\int_0^t C(\vv) (1+|x|^2) \e^{ \ll_1 s}  \e^{\ll_\vv(t-s)}\d s\\
&\le  \ff{\ll_1 C(\vv)}{\ll_1-\ll_\vv}  \e^{\ll_1 t} (1+|x|^2),\ \ x\in\H,\ \ t\ge 0.\end{align*} Therefore, the proof is finished.
\end{proof}

\beg{lem} \label{L2}  Under $(H_1)$ and $(H_2)$, there exists a constant $k>0$ such that
\beq\label{NB} \int_\H \ff{\mu(\d x)}{\int_\H \e^{-\ll |x-y|^2 }\mu(\d y)}\le   \e^{\gg(k\ll)}
 \prod_{i=1}^\infty \ff{\ll_i+k\ll}{\ss{\ll_i^2-\ff 1 2\ll_1^2}},\ \ \ll\ge 1.\end{equation}
 \end{lem}

\beg{proof}
Let $\{e_i\}_{i\ge 1}$ be the eigen-basis of
$A$, i.e. it is an orthonormal basis of   $\H$ such that
$$A e_i=\ll_i e_i,\ \ i\ge 1.$$ Each $x\in\H$ is corresponding to an eigen-coordinate
$$ (x_i)_{i\ge 1}:= (\<x, e_i\>)_{i\ge 1}\in \ell^2:=\Big\{(r_i)_{i\ge 1}\subset \R^\infty: \sum_{i=1}^\infty r_i^2<\infty\Big\}.$$
Under this coordinate we have 
\beq\label{*DP}  \mu_0(\d x)= \prod_{i=1}^\infty \ff{\ss{\ll_i}}{\ss{2\pi}} \e^{-\ff{\ll_i x_i^2} 2}\d x_i. \end{equation}
Combining this with $(H_2)$ and $\mu(\d x)=Z_V^{-1}\e^{V(x)}\mu_0(\d x)$,
 we obtain
 \beg{align*} &I:= \int_\H \ff{\mu(\d x)}{\int_\H \e^{-\ll |x-y|^2 }\mu(\d y)}
\le  \e^{\gg(\vv^{-1})}
 \prod_{i=1}^\infty \bigg\{\int_\R \ff{\e^{-\ff{\ll_i-\vv}2 x_i^2}} {\int_\R \e^{-\ll |x_i-y_i|^2 -\ff{\ll_i+\vv}2 y_i^2}\d y_i}\bigg\}\d x_i,\ \ \vv>0.\end{align*} 
 Noting that 
 $$\ll|x_i-y_i|^2 +\ff{\ll_i+\vv} 2 y_i^2 = \ff{2\ll +\ll_i+\vv} 2 \Big(y_i-\ff{2\ll x_i}{2\ll+\ll_i +\vv}\Big)^2 +\ff{\ll(\ll_i +  \vv)x_i^2}{2\ll+\ll_i+\vv},$$ we derive
\beq\label{*YPP} \int_\R \e^{-\ll|x_i-y_i|^2 -\ff{\ll_i+\vv}2 x_i^2} \d y_i =\Big(\ff {2 \pi}{2\ll+\ll_i +\vv}\Big)^{\ff 1 2} \e^{\ff{-\ll(\ll_i +  \vv)x_i^2}{2\ll+\ll_i+\vv}},\ \ \vv> 0.\end{equation}
So, for $\vv>0$ such that $\ll_1^2-4\ll\vv-\vv^2>0$, we have 
\beq\label{**P} \beg{split} 
 I &\le \e^{\gg(\vv^{-1})}\prod_{i=1}^\infty\bigg\{\Big(\ff{2\ll + \ll_i+\vv}{4\pi}\Big)^{\ff 1 2} \int_\R \e^{-(\ff{\ll_i-\vv }2- \ff{\ll\ll_i +\ll\vv}{2\ll+\ll_i+\vv}) x_i^2} \d x_i\bigg\}\\
 &=\e^{\gg(\vv^{-1})}\prod_{i=1}^\infty\bigg\{\Big(\ff{2\ll + \ll_i+\vv}{4\pi}\Big)^{\ff 1 2} \Big(\ff{4\pi(2\ll+\ll_i+\vv)}{\ll_i^2-4\ll\vv-\vv^2}\Big)^{\ff 1 2}\bigg\}, \ \ \vv >0. \end{split} \end{equation} 
Taking $\vv:= \ss{4\ll^2+\ff 1 2\ll_1^2}-2\ll,$ we have 
\beg{align*} &\vv\in \bigg(\ff{\ll_1^2}{2\ss{16\ll^2+2\ll_1^2}}, \ \ff{\ll_1^2}{8\ll}\bigg),\\ 
&\ll_i^2-4\ll\vv-\vv^2\ge \ll_i^2-\ff 1 2 \ll_1^2>0,\ \ i\ge 1.\end{align*} 
Then  for any $\ll\ge 1$, 
$$\vv^{-1}\le \ff{2\ss{16+2\ll_1^2}}{\ll_1^2}  \ll,\ \ \ \ \Big(\ff{2\ll + \ll_i+\vv}{4\pi}\Big)^{\ff 1 2} \Big(\ff{4\pi(2\ll+\ll_i+\vv)}{\ll_i^2-4\ll\vv-\vv^2}\Big)^{\ff 1 2}\le  \ff{\ll_i+(2+\ff{\ll_1^2}8)\ll}{\ss{\ll_i^2-\ff 1 2\ll_1^2}}.$$
Therefore, \eqref{**P} implies \eqref{NB} for some constant $k>0$.    
\end{proof}

\beg{proof}[Proof of Corollary \ref{C1.2}]  (1) By   \eqref{A0} and the second formula in \eqref{AA}, we find a constant $c_1>1$ such that
$$\bb(\vv)\le \int_{\H} \ff{\mu(\d x)\mu(\d y)}{\int_\H\e^{-c_1\vv^{-1} |x-y|^2}},\ \ \vv\in (0,1).$$
Combining this with   Lemma \ref{L2}, we find constants $c_2,c_3, c_4>0$ such that
\beg{align*} \bb(\vv)&\le \e^{\gg( c_2\vv^{-1})} \exp\bigg[\sum_{i=1}^\infty \log\bigg(1+ \ff{\ll_i+c_2\vv^{-1}- (\ll_i^2-\ff 1 2\ll_1^2)^{\ff 1 2}}{(\ll_i^2-\ff 1 2 \ll_1^2)^{\ff 1 2}}\bigg)\bigg]\\
&\le \e^{\gg( c_2\vv^{-1})} \exp\bigg[c_3\vv^{-1} \sum_{i=1}^\infty \ff 1 {\ll_i}\bigg]\le \e^{\gg(c_2\vv^{-1})+ c_4\vv^{-1}},\ \ \vv\in (0, c_1).\end{align*}  Noting that $\bb(\vv)$ is decreasing in $\vv$, we find a constant $k>0$ such that
\beq\label{BB0} \bb(\vv)\le \e^{\gg(k\vv^{-1})+ k\vv^{-1}},\ \ \vv\in (0,1).\end{equation}
On the other hand, by the definition of the mild solution and that of $\aa$ in \eqref{AA}, we have
\beq\label{BB2} \beg{split} &\aa(\vv)=\E^\mu\big[ |X_\vv-X_0|^2\big] \\
&= \E^\mu\bigg[\bigg|\e^{-A\vv}X_0-X_0+ \int_0^\vv \e^{-A(\vv-s)}\nn V(X_s)\d s + \ss 2 \int_0^\vv \e^{-A(\vv-s)}\d W_s\bigg|^2\bigg],\\
&\le 3 \E^\mu \big[|\e^{-A\vv}X_0-X_0|^2\big]+3\vv\int_0^\vv \E^\mu\big[|\nn V(X_s)|^2\big]\d s +
6\int_0^\vv \|\e^{-A(\vv-s)}\|_{HS}^2\d s.
\end{split} \end{equation}
Moreover, by  $(H_1)$ and $(H_2)$,  $\nn V(x)$ is Lipschitz continuous hence has  a linear growth in $|x|$,  and   $\mu(|\cdot|^2)<\infty.$ So, \eqref{A00} implies
  $\sup_{s\in [0,1]} \E^\mu[|\nn V(X_s)|^2]<\infty$. Thus, by $(H_1)$ and $(H_2)$ which imply
  $$\E^\mu[\<X_0,e_i\>^2]=\mu(|x_i|^2)\le \ff c {\ll_i},\ \ i\ge 1$$ for some constant $c>0$,    we find   constants $c_5,c_6>0$ such that
\beg{align*} &\vv \int_0^\vv \E [| \nn V(X_s)|^2]\d s  +  \E^\mu[|\e^{-A\vv}X_0-X_0|^2] + \int_0^\vv \|\e^{-A(\vv-s)}\|_{HS}^2\d s\\
&=c_5\vv^2 + \sum_{i=1}^\infty \bigg(\ff{(1-\e^{-\ll_i\vv})^2}{\ll_i} + \int_0^\vv \e^{-2\ll_i(\vv-s)}\d s\bigg)
\le c_6\sum_{i=1}^\infty \ff{1-\e^{-2\ll_i\vv}}{\ll_i},\ \ \vv\in (0,1).\end{align*}
Substituting into \eqref{BB2},  we find  a constant    $k>0$ such that
\beq\label{AAP} \aa(\vv)\le k \sum_{i=1}^\infty \ff{1-\e^{-2\ll_i\vv}}{\ll_i},\ \ \vv\in (0,1).\end{equation}
Combining this with \eqref{BB0} and applying Theorem \ref{T1.1}, we prove (1).

(2) According to the first assertion and \eqref{A2}, it suffices to show that for any $k>K^+$ there exist  constants $r, k_r>0$ such that
$$ c(r,x)\le k_r \e^{k_r|x|^2},\ \ x\in \H,$$ which follows from $(H_2)$ and  \eqref{*YPP} with $\ll=\ff{2K}{1-\e^{-2Kr}} $ and $\vv=0$. 

\end{proof}

\paragraph{Example 2.1.}  Let $\nn V$ be Lipschitz continuous, $\ll_i\ge c_0 i^p$ for some constants $c_0>0$ and $p>1$,   and there exists a constant $c>0$   such that
\beq\label{V2'} |V(x)|\le c(1+|x|),\ \ x\in \H\end{equation} holds. Then  there exists a a constant $\kk>0$ such that
\beq\label{FF1} \E^\mu [\W_2(\mu_t,\mu)^2] \le \kk (\log t)^{p^{-1}-1},\ \ t\ge 2.\end{equation}
If moreover $(H_3)$ holds, then for any $k>K^+$ there exists a constant $c(k)>0$ such that
\beq\label{FF2} \big(\E^\mu [\W_2(\mu_t,\mu)] \big)^2\le c(k) \e^{k|x|^2}  (\log t)^{p^{-1}-1},\ \ t\ge 2,\ \ x\in \H.\end{equation}

\beg{proof} Let
\beq\label{BCC} h(\vv)= \sum_{i=1}^\infty\ff{1-\e^{-2\vv \ll_i}}{\ll_i},\ \ \vv\in [0,1].\end{equation}  When $\ll_i\ge c i^p$ for some constants $c>0$ and $p>1$,
we find a constant $c_1>0$ such that
$$ h'(\vv)= \sum_{i=1}^\infty 2\e^{-2\vv \ll_i}\le 2 +2\int_1^\infty \e^{-2c\vv s^p}\d s \le c_1\vv^{-  p^{-1}},\ \ \vv \in (0,1].$$
Thus, there exists a constant $c_2>0$ such that
\beq\label{BCC0} \sum_{i=1}^\infty\ff{1-\e^{-2\vv \ll_i}}{\ll_i}= \int_0^\vv h'(s)\d s\le c_2 \vv^{1-p^{-1}},\ \ \vv\in (0,1].\end{equation}
On the other hand,  \eqref{V2'} implies $(H_2)$ with
$$\gg(s)= c_3s,\ \ s\ge 1$$ for some constant $c_3>0$. Then  by taking $\vv=\ff{2(c_3+k)}{\log t}$,  we find constants $c_4,c_5>0$ such that
\beg{align*} &\inf_{\vv\in (0,1)}   \bigg\{ \ff 1 t \e^{k\vv^{-1}+\gg(k\vv^{-1})}  +\sum_{i=1}^\infty
\ff{1-\e^{-2\ll_i\vv}}{\ll_i} \bigg\}\\
&\le c_4 \inf_{\vv\in (0,1)} \Big\{\ff 1 t \e^{(k+c_3)\vv^{-1}}+ c_2 \vv^{1- p^{-1}}\Big\}\\
&\le c_5 (\log t)^{p^{-1}-1},\ \ t\ge 2.\end{align*}
Therefore, the desired assertions follow from Corollary \ref{C1.2}.
\end{proof}

\paragraph{Example 2.2.} Let $\nn V$ be Lipschitz continuous,  $\ll_i\ge c \e^{i^p }$ for some constant $c>0$ and $p>0$, and  \eqref{V2'} holds for some constant $c>0$. Then
  there exists a a constant $\kk>0$ such that
\beq\label{FF1'} \E^\mu [\W_2(\mu_t,\mu)^2] \le \kk (\log t)^{-1} \log\log t,\ \ t\ge 4.\end{equation}
If moreover $(H_3)$ holds, then for any $k>K^+$ there exists a constant $c(k)>0$ such that
\beq\label{FF2'} \E^\mu [\W_2(\mu_t,\mu)^2] \le c(k) \e^{k|x|^2}  (\log t)^{-1}\log\log t,\ \ t\ge 4,\ \ x\in \H.\end{equation}

\beg{proof} Let $h$ be in \eqref{BCC}. When $\ll_i\ge c  \e^{c i^p }$ for some constant $c>0$ and $p>0$, by using the integral transform $r=\vv \e^{cs^p}$, we find constants $c_1,c_2>0$ such that
\beg{align*}& h'(\vv)= 2\sum_{i=1}^\infty 2\e^{-2\vv \ll_i}\le  2\int_0^\infty \exp\big[-2c\vv \e^{cs^p}\big]\d s \\
&= 2\int_{\vv}^\infty  \e^{-2cr} \ff{\d }{\d r} \big\{c^{-1}\log[r \vv^{-1}]\big\}^{\ff 1 p}\,\d r \\
&\le   c_1 \int_\vv^1   \big\{\log r+\log \vv^{-1}\big\}^{\ff 1 p -1} \d\log r +
c_1 \big\{\log (1+\vv^{-1})\big\}^{\ff 1 p -1}\\
&=c_1 \log (1+\vv^{-1})^{\ff 1 p-1} + c_0 \int_{\log \vv}^0   \big\{u+\log \vv^{-1}\big\}^{\ff 1 p -1} \d  u\\
&\le c_2 \log (1+\vv^{-1})^{\ff 1 p},\ \ \vv \in (0,1].\end{align*}
Thus, there exists a constant $c_3>0$ such that
$$\sum_{i=1}^\infty\ff{1-\e^{-2\vv \ll_i}}{\ll_i}= \int_0^\vv h'(s)\d s\le c_3 \vv  \log (1+\vv^{-1})^{\ff 1 p},  \ \ \vv\in (0,1].$$  So, as in the proof of Example 2.1,  we find constants
$c_4, c_5>0$ such that
\beg{align*} &\inf_{\vv\in (0,1)}   \bigg\{ \ff 1 t \e^{k\vv^{-1}+\gg(k\vv^{-1})}  +\sum_{i=1}^\infty
\ff{1-\e^{-2\ll_i\vv}}{\ll_i} \bigg\}\\
&\le c_4 \inf_{\vv\in (0,1)} \Big\{\ff 1 t \e^{(k+c_4)\vv^{-1}}+  c_3 \vv  \log (1+\vv^{-1})^{\ff 1 p} \Big\}\\
&\le c_5 (\log t)^{-1}(\log\log t)^{\ff 1 p} ,\ \ t\ge 4.\end{align*}
Therefore, the desired assertions follow from Corollary \ref{C1.2}.
\end{proof}

\section{Lower bound estimate}

We first present a lower bound estimate on 
\beq\label{LB0} \W_p(\mu,\nu):= \inf_{\pi\in \C(\mu,\nu)} \bigg\{\int_{E\times E} \rr(x,y)^p \pi(\d x,\d y)\bigg\}^{\ff 1 p},\ \
p>0, \mu,\nu\in \scr P(E)\end{equation} 
 for a metric space  $(E,\rr)$, where  $\scr P(E)$ is the set of  all probability measures on $E$.
 As a generalization to \cite[Proposition 4.2]{RE1} which essentially works for the finite-dimensional setting, we have the following result which also applies to   infinite dimensions.

\beg{lem}\label{LB1} Let $\mu\in \scr P(E) $ such that
\beq\label{LB1} \sup_{x\in E}\mu(B(x,r)) \le \psi(r),\ \ r\ge 0 \end{equation} holds for some increasing function $\psi,$ where
$B(x,r):=\{y\in E: \rr(x,y)<r\}$. Then
for any $N\ge 1$ and any probability measure $\mu_N$ supported on a set of $N$ points in $E$,
\beq\label{LB2} \W_p(\mu_N,\mu) \ge 2^{-\ff 1 p}  \psi^{-1}\Big(\ff 1 {2N}\Big),\end{equation}
where $\psi^{-1}(s):=\sup\{r\ge 0: \psi(r)\le s\}, s\ge 0.$
\end{lem}
\beg{proof} Let $D={\rm supp}\mu_N$ which contains $N$ many points, so that from \eqref{LB1} we conclude that
 $D_r:= \cup_{x\in D} B(x,r)$ satisfies
$$\mu(D_r)\le \sum_{x\in D} \mu(B(x,r))\le N \psi(r),\  \ r\ge 0.$$
Therefore, for any $\pi\in \C(\mu_N, \mu)$, we get
$$ \int_{E\times E} \rr(x,y)^p \pi(\d x,\d y) \ge \int_{D\times D_r^c} r^p \pi(\d x,\d y) = r^p \mu(D_r^c) \ge r^p\{1-N\psi(r)\},\  \ r\ge 0.$$
Combining this with \eqref{LB0} we obtain
$$\W_p(\mu,\nu)^p \ge \sup_{r\ge 0} r^p[1-N\psi(r)]\ge \ff 1 2 \big\{\psi^{-1}(1/(2N))\big\}^p.$$
\end{proof}

Now, let $E=\H$ and consider 
$$\tt\W_1(\mu,\nu)= \inf_{\pi\in \C(\mu,\nu)} \int_{\H\times \H} \{|x-y|\land 1\} \pi(\d x,\d y),\ \ \mu,\nu\in \scr P.$$

\beg{thm}\label{T2.1} Assume $(H_1)$. Then there exists a constant $k>0$ such that
\beq\label{LB3} \E^\mu[\tt\W_1(\mu_t,\mu)] \ge \sup_{N\in \mathbb N} \bigg\{\ff 1 2 \psi^{-1}\big((2N)^{-1}\big) -\ff k N - k \Big(\sum_{i=1}^\infty \ff{1-\e^{-2\ll_i t/N}}{\ll_i}\Big)^{\ff 1 2}\bigg\},\ \ t\ge 1.\end{equation}
If moreover $(H_2)$ holds, then there exists a constant $k>0$ such that for any $x\in \H$,
\beq\label{NMM0}  \E^x[\tt\W_1(\mu_t,\mu)] 
 \ge \sup_{N\in \mathbb N} \bigg\{\ff 1 2 \psi^{-1}\big((2N)^{-1}\big)  - \Big(k (1+|x|^2)\sum_{i=1}^\infty \ff{1-\e^{-2\ll_i t/N} }{\ll_i}\Big)^{\ff 1 2}\bigg\},\ \ t>0.  \end{equation}
\end{thm}

\beg{proof} For any $t>0$ and $N\in \mathbb N$, let
$$t_i= \ff{(i-1)t}N,\ \ 1\le i\le N+1.$$
Take
$$\mu_{t,N} = \ff 1 N \sum_{i=1}^N\dd_{X_{t_i}} =\ff 1 t \sum_{i=1}^N \int_{t_i}^{t_{i+1}} \dd_{X_s} \d s.$$
Noting that
$$\pi(\d x,\d y) := \ff 1 t \sum_{i=1}^N \int_{t_i}^{t_{i+1}} \dd_{X_{t_i}}(\d x) \dd_{X_{t_i}}(\d y) \d s\in \C(\mu_{t,N},\mu_t),$$
we obtain
\beq\label{LBBN} \beg{split}  &\E^\mu \big[\tt\W_1(\mu_{t,N},\mu_t)\big] \le \E^\mu \bigg[\ff 1 t \sum_{i=1}^N \int_{t_i}^{t_{i+1}}\big\{|X_s-X_{t_i}|\land 1\big\}\, \d s\bigg]\\
&\le \ff 1 t \sum_{i=1}^N \int_{t_i}^{t_{i+1}}\E^\mu [|X_0-X_{s-t_i}|] \d s\le \sup_{s\in [0, t/N]} \big(\E^\mu[|X_0-X_s|^2]\big)^{\ff 1 2}.\end{split}\end{equation}
This together with \eqref{AAP} implies
$$\E^\mu \big[\tt\W_1(\mu_{t,N},\mu_t)\big]\le \bigg(k\sum_{i=1}^\infty\ff{1-\e^{-2\ll_i t/N}}{\ll_i}\bigg)^{\ff 1 2}.$$
Therefore, by combining with \eqref{LB2} for $E=\H$ and $\rr(x,y)=|x-y|\land 1$, we arrive at
\beq\label{NMM} \beg{split} &\E^\mu \big[\tt\W_1(\mu,\mu_t)\big]\ge \E^\mu\big[\tt\W_1(\mu_{t,N},\mu)\big]- \E^\mu \big[\tt\W_1(\mu_{t,N},\mu_t)\big]\\
&\ge \ff 1 2 \psi^{-1}\Big(\ff 1 {2N}\Big) - \bigg(k\sum_{i=1}^\infty\ff{1-\e^{-2\ll_i t/N}}{\ll_i}\bigg)^{\ff 1 2},\ \ N\in\mathbb N.\end{split}
\end{equation}
Then \eqref{LB3} holds.

On the other hand,  by  Lemma \ref{L1} and the linear   growth of $|\nn V|$,  we find   constants  $c_1,c_2 >0$ such that
for any $s\in [t_i,t_{i+1}]$,
\beq\label{CW1}  \beg{split} &\E^x |X_s-X_{t_i}| - \E^x\big|  (\e^{-(s-t_i)A}-1)X_{t_i} \big|\\
&\le  \E^x\bigg|   \int_{t_i}^s \e^{-A(r-t_i) }\nn V(X_r)\d r +\ss 2 \int_{t_i}^s \e^{-A(r-t_i)} \d W_r\bigg|\\
&\le   \int_{t_i}^s \e^{-\ll_1(r-t_{i})} \E^x[|\nn V(X_r)|]\d r+ \ss 2 \bigg(\sum_{j=1}^\infty \int_{t_i}^s \e^{-2\ll_j (r-t_i)}\d r\bigg)^{\ff 1 2}\\
&\le   \ff{c_1 (1-\e^{-\ll_1 N/t})}{\ll_1}  (1+|x|) + \ss 2 \bigg(\sum_{j=1}^\infty \ff{1-\e^{-2\ll_j t/N}}{2\ll_j}\bigg)^{\ff 1 2}\\
&\le   c_2 (1+|x|)  \bigg(\sum_{j=1}^\infty \ff{1-\e^{-2\ll_j t/N}}{ \ll_j}\bigg)^{\ff 1 2}.\end{split} \end{equation} where the last step follows from the fact that 
$$\bigg(\sum_{j=1}^\infty \ff{1-\e^{-2\ll_j t/N}}{ \ll_j}\bigg)^{\ff 1 2}\ge \bigg(  \ff{1-\e^{-2\ll_1 t/N}}{ \ll_1}\bigg)^{\ff 1 2} \ge  \ff{1-\e^{-2\ll_1 t/N}}{\ss{ \ll_1}}.$$
Moreover, by $s-t_i\in [0, t/N]$ and Jensen's inequality, 
\beq\label{CW2} \beg{split} & \E^x\big|  (\e^{-(s-t_i)A}-1)X_{t_i} \big| \le \E^x \ss{\sum_{j=1}^\infty (1- \e^{-t\ll_j/N})^2\<X_{t_i},e_j\>^2}\\
&\le \bigg(\sum_{j=1}^\infty (1- \e^{-t\ll_j/N})^2\E^x \<X_{t_i},e_j\>^2\bigg)^{\ff 1 2}.\end{split}\end{equation} 
Since for the  $1$-dimensional Brownian  motion $B_s^j:=\<W_s, e_j\>$ we have 
\beg{align*} \d \<X_s,,e_j\>^2 &= \big(-2\ll_j \<X_s,e_j\>^2+2\<\nn V(X_s), e_j\> \<X_s,e_j\>\big)\d s+ 2 \ss 2 \<X_s,e_j\> \d B_s^j\\
&\le \big(- \ll_j \<X_s,e_j\>^2+\ll_j^{-1} |\nn V(X_s)|^2 \big)\d s  + 2 \ss 2 \<X_s,e_j\> \d B_s^j,\end{align*} 
by  Lemma \ref{L1} and the linear   growth of $|\nn V|$,  we find a constant   $c _3 >0$ such that
$$\E^x \<X_{t_i},e_j\>^2 \le \<x,e_j\>^2 \e^{-\ll_j t_i} + \ff {c_3}  {\ll_j^2} (1+|x|^2),\ \ j\ge 1.$$ 
Combining this with \eqref{CW1} and \eqref{CW2} we find a constant  $c_4>0$ such that
\beq\label{R1}  \E^x |X_s-X_{t_i}| \le  c_4 (1+|x|)    \bigg(  \sum_{j=1}^\infty\Big((1-\e^{-\ll_j t/N})^2 \e^{-\ll_j t_i}+ \ff{1-\e^{-2\ll_j t/N}}{ \ll_j}\Big)\bigg)^{\ff 1 2}.\end{equation} 
By the same reason leading to \eqref{LBBN} we have 
$$\W_1(\mu_{N,t},\mu_t) \le \ff 1 t  \sum_{i=1}^N \int_{t_i}^{t_{i+1}}|X_s-X_{t_i}| \d s\le \ff 1 N \int_0^{t/N} |X_s-x|\d s + \ff 1 t  \sum_{i=2}^N \int_{t_i}^{t_{i+1}}|X_s-X_{t_i}| \d s.$$
This together with  Lemma \ref{L1} and \eqref{R1} yields that for some $c_5 >0$ and  any $t\ge 1, N\in \mathbb N,$
\beg{align*} &\E^x[\W_1(\mu_{t,N},\mu_t)] \le \ff{c_5}N (1+|x|) + \ff 1 t\sum_{i=2}^N  \int_{t_i}^{t_{i+1}} \E^x[|X_s-X_{t_i}|]\d s\\
&\le  c_5 (1+|x|) \bigg[\ff 1 N +\ff 1 N \sum_{i=2}^N   \Big(\sum_{j=1}^\infty \e^{-\ll_j (i-1)t/N}\big(1-\e^{-\ll_j t/N}\big)^2\Big)^{\ff 1 2} +  \bigg(\sum_{j=1}^\infty \ff{1-\e^{-2\ll_jt/N}}{ \ll_j}\bigg)^{\ff 1 2}\bigg]\\
&\le c_5 (1+|x|) \bigg[\ff 1 N +  \Big(\ff 1 N \sum_{i=2}^N  \sum_{j=1}^\infty \e^{-\ll_j (i-1)t/N}\big(1-\e^{-\ll_j t/N}\big)^2\Big)^{\ff 1 2} +  \bigg(\sum_{j=1}^\infty \ff{1-\e^{-2\ll_jt/N}}{ \ll_j}\bigg)^{\ff 1 2}\bigg]\\
&\le  c_5 (1+|x|) \bigg[ \ff 1 N+   2 \bigg(\sum_{j=1}^\infty \ff{1-\e^{-2\ll_jt/N}}{ \ll_j}\bigg)^{\ff 1 2}\bigg],  
 \end{align*} where the last step follows from
 $$\ff 1 N \sum_{i=2}^N    \e^{-\ll_j (i-1)t/N} \le \ff 1 t \int_0^t \e^{-\ll_j s}\d s \le \ff 1 {\ll_j}.$$
Therefore, as in \eqref{NMM} we prove \eqref{NMM0} for some constant $k>0$.

\end{proof}

\paragraph{Example 3.1.} Assume $(H_1), (H_2)$.    If there exist constants $  p\ge q>1$ and $k_1,k_2>0$ such that
\beq\label{XZ0} k_1i^q\le \ll_i\le k_2 i^p,\ \ i\ge 1,\end{equation}
then there exists a constant $c>0$ such that for large $t>1$,
\beq\label{XZ1} \E^\mu [\tt\W_1(\mu_t,\mu)] \ge c
\{\log t\}^{-(\ff{p-1}2\land 1)}.  \end{equation}
Moreover,  for any $x\in \H$ there exist constants $c(x),t(x)>0$ such that  
\beq\label{XZ2} \E^x [\tt\W_1(\mu_t,\mu)] \ge c(x)
\{\log t\}^{-\ff{p-1}2\land 1}, \ \ t\ge t(x).  \end{equation}

\beg{proof} (a) We first consider $p>2$. By $(H_2)$ with $\vv=\ff{\ll_1}2$ and \eqref{*DP},  we find a constant $c_1>0$ such that
\beg{align*} &\psi(r):= \sup_{x\in\H}\mu(B(x,r)) \le    c_1 \prod_{i=1}^\infty \ff{\ss{\ll_i}}{\ss{2\pi}} \int_0^r \e^{-\ff{(\ll_i-\vv) s^2}2}\d s\\
&=  c_1 \prod_{i=1}^\infty \ff{\ss{\ll_i}}{\ss{\ll_i-\vv}} \bigg(1- \ff{\ss{ \ll_i-\vv }}{\ss{2\pi}}\int_r^\infty \e^{-\ff{(\ll_i -\vv)s^2}2}\d s\bigg).\end{align*}
Since $\ll_i\ge k_2i^p$ and  $p>2$,  and $\vv=\ff {\ll_1}2\le \ff{\ll_i}2$, there exist   constants $c_2, c_3, c_4>0$ such that 
\beg{align*} &\log \Big[\ff{\psi(r)}{c_1}\Big]\le  \sum_{i=1}^\infty\bigg\{ \ff 1 2 \log \Big(1+\ff{ \vv }{ \ll_i-\vv } \Big)+\log\bigg( 1-\ff{\ss{\ll_i-\vv}}{\ss{2\pi}}  \int_r^\infty 
\e^{-\ff{(\ll_i -\vv)s^2}2}\d s\bigg)\bigg\}\\
&\le      \sum_{i=1}^\infty \bigg\{\ff{\vv}{2(\ll_i-\vv)} - \ff{\ss{\ll_i-\vv}}{\ss{2\pi}} \int_r^\infty 
\e^{-\ff{(\ll_i -\vv)s^2}2}\d s \bigg\}\\
&  \le      \sum_{i=1}^\infty\Big( \ff 1  { \ll_i } -      \ff{c_2}{\ss{\ll_i}}  \Big)\le c_3 - c_4 r^{-1},\ \ \ r>0. 
\end{align*}   Since $\psi(r)\le 1$, this implies that for some constant $c_5>0$, 
$$\psi(r)\le c_5\e^{- c_4 r^{-1}},\ \ r>0.$$ 
 Therefore, there exist     constants  $c_6>0$    such that
\beq\label{GPP} \psi^{-1}(1/(2N)) \ge c_6 \{\log (2N)\}^{-1},\ \ N\ge 1. \end{equation}

  On the other hand, since $\ll_i\ge k_1 i^q$ for some $q>1$,  \eqref{BCC0} holds for $q$ replacing $p$, i.e. there exists a constant $k>0$ such that
  \beq \label{BCC0'} \sum_{i=1}^\infty\ff{1-\e^{-2\vv \ll_i}}{\ll_i} \le k \vv^{1-q^{-1}},\ \ \vv\in (0,1]. \end{equation}
   Combining this with \eqref{GPP} and   \eqref{LB3} with   $N= 1+ \lceil t\rceil^2$,   where $  \lceil t\rceil$ is the integer part of $t$, we find a constant $c>0$ such that for large $t$
  \beq\label{XZ1-1} \E^\mu [\tt\W_1(\mu_t,\mu)]  \ge  c
\{\log t\}^{-1}. \end{equation}
Similarly, for any $x\in\H$ there exist constants $c(x),t(x)>0$ such that
 \beq\label{XZ2-1} \E^x [\tt\W_1(\mu_t,\mu)]  \ge  c_1
\{\log t\}^{-1},\ \ t\ge t(x). \end{equation}

(b) Take $\vv=1$ in $(H_2)$, we find a constant $C_1>0$ such that for any $R>0$,  
\beg{align*} &\psi(r)\le C_1 \int_{B(0,r)} \e^{\ff{|x|^2}2} \mu_0(\d x) \\
&\le C_1 \e^{\ff{(R+1)r^2}2} \int_\H \e^{-\ff R 2 |x|^2}\mu_0(\d x) = C_1 \e^{\ff{(R+1)r^2}2} \prod_{i=1}^\infty
\ff{\ss{\ll_i}}{\ss{\ll_i+R}}.\end{align*}  Since $\ll_i\le k_2i^p$ for $i\ge 1$, this implies
 \beg{align*} &\log\Big[\ff{\psi(r)}{C_1}\Big]\le \ff{(R+1)r^2}2 +\ff 1 2 \log\Big(1-\ff{R}{\ll_i+R}\Big)\le \ff{Rr^2}2 -\ff 1 2 \sum_{i=1}^\infty \ff{R}{\ll_i+R}\\
&\le \ff{(R+1)r^2}2 -\ff R 2 \int_1^\infty \ff {\d s}{k_2 s^p+R} \le \ff{(R+1)r^2}2 -c_1 R(1+R^{\ff 1 p})^{1-p} \\
&\le  Rr^2  - c_2 R^{\ff 1 p},\ \ r>0, R\ge 1\end{align*}
for some constants $c_1,c_2>0$. By taking $r_0:= \vv^{\ff{p-1}{2p}}$ and $R= \vv r^{-\ff{2p}{p-1}}$ for small enough $\vv>0$, we find a constant $c_3>0$ such that  
$$\log\Big[\ff{\psi(r)}{C_1}\Big]\le Rr^2 - c_2 R^{\ff 1 p}\le -c_3 r^{-\ff 2 {p-1}},\ \ r\in (0,r_0].$$
Combining this with $\psi(r)\le 1$ for all $r\ge 0$,  we find a constant $c_4>0$ such that
$$\psi(r)\le c_4 \e^{-c_3 r^{-\ff 2 {p-1}}},\ \ r>0.$$ This implies
$$\psi^{-1}(1/(2N)) \ge c_5 \{\log (2N)\}^{-\ff{p-1}2},\ \ N\in \mathbb N$$
for some constant $c_5>0$. Combining this with \eqref{LB3}, \eqref{BCC0'} and taking  $N= 1+ \lceil t\rceil^2$ for large $t>0$, we find a constant $c_6>0$ such that
$$\E[\tt\W_1(\mu_t,\mu)]\ge c_6 \{\log t\}^{-\ff {p-1}2} $$ holds for large $t>0$. 
This together with \eqref{XZ1-1} implies \eqref{XZ1}. Similarly, \eqref{XZ2} holds for any $x\in\H$ and some constants $c(x), t(x)>0.$
 \end{proof}

\paragraph{Acknowledgement.}  The author  would like to thank   the referee  for helpful comments and corrections.

  \end{document}